\documentclass{amsart}

\usepackage{amsmath,amssymb,amsthm}

\theoremstyle{definition}

\newtheorem*{ack}{Acknowledgements}

\def\co{\colon\thinspace}

\newcommand{\C}{\mathbb{C}}
\newcommand{\R}{\mathbb{R}}

\newcommand{\rmd}{\mathrm{d}}
\newcommand{\rme}{\mathrm{e}}
\newcommand{\bfx}{\mathbf{x}}
\newcommand{\bfy}{\mathbf{y}}

\newcommand{\alst}{\alpha_{\mathrm{st}}}
\newcommand{\omst}{\omega_{\mathrm{st}}}

\begin{document}

\author[H.~Geiges]{Hansj\"org Geiges}
\address{Mathematisches Institut, Universit\"at zu K\"oln,
Weyertal 86--90, 50931 K\"oln, Germany}
\email{geiges@math.uni-koeln.de}

\author[N.~R\"ottgen]{Nena R\"ottgen}
\author[K.~Zehmisch]{Kai Zehmisch}
\address{Mathematisches Institut, WWU M\"unster,
Einstein\-stra\-{\ss}e 62, 48149 M\"unster, Germany}
\email{roettgen@uni-muenster.de}
\email{kai.zehmisch@uni-muenster.de}

\title[A Hamiltonian plug]{From a Reeb orbit trap to a Hamiltonian plug}

\date{}

\begin{abstract}
We present a simple construction of a plug for Hamiltonian
flows on hypersurfaces of dimension at least five
by doubling a trap for Reeb orbits.
\end{abstract}


\maketitle


\section{Introduction}
In \cite{seif50}, H.~Seifert remarked that `it is unknown if every
continuous [nowhere vanishing] vector field of the three-dimensional
sphere $S^3$ contains a closed integral curve.' The supposition that
the answer to this question is positive has become known as the
Seifert conjecture. In its original form, this conjecture
was disproved by P.~A.~Schweitzer~\cite{schw74}. Earlier,
F.~W.~Wilson~\cite{wils66} had shown the existence of aperiodic flows
on any compact manifold of dimension at least four
and of vanishing Euler characteristic (which
is the condition for the existence of a non-singular vector field).
His result is based on the construction of what,
following J.~Harrison~\cite{harr88} and K.~Kuperberg~\cite{kupe-k94},
is now known as a \emph{plug}. This is a local model of an aperiodic flow
that can be inserted into a flow box around an isolated periodic orbit
of a given flow in order to open up that orbit. Needless to say,
care has to be taken not to create new periodic orbits
by this process. For a beautiful survey on constructions of
aperiodic flows see~\cite{kupe-k99}.

In this note, we are concerned with the Hamiltonian version of the Seifert
conjecture: does the Hamiltonian flow on a closed hypersurface
in a symplectic manifold necessarily have a periodic orbit?
In this generality, the conjecture has been disproved
for hypersurfaces of dimension at least five by
V.~Ginzburg~\cite{ginz95,ginz97} and M.~Her\-man~\cite{herm99,herm};
Ginzburg's construction was simplified by E.~Kerman~\cite{kerm02}.
In dimension three, the best counterexample to date is one where
the Hamiltonian function defining the hypersurface
is only $C^2$-smooth, and hence the Hamiltonian vector field only
$C^1$-smooth~\cite{gigu03}.
By contrast, there are many positive results
under more restrictive assumptions, such as the hypersurface being of
contact type, or for dense subsets of levels of a Hamiltonian function.
These matters are surveyed comprehensively in~\cite{ginz99};
see also \cite{herm} and the introductions to \cite{gigu03,kerm02}.

Here we describe a simple construction of a plug for smooth Hamiltonian
flows on hypersurfaces of dimension at least five, starting from
a trap for Reeb orbits invented in~\cite{grz14}. 

\section{Plugs}
For our purposes, a \emph{plug} will be a certain smooth non-singular
vector field $X$ on $D\times I$,
where $D=D^{m-1}_{\delta}$ is an $(m-1)$-dimensional disc of
radius~$\delta$, and $I$ an interval $[-\varepsilon,\varepsilon]$.
Write $z$ for the $I$-coordinate.
This vector field is supposed to have the following properties:
\begin{itemize}
\item[(i)] $X=\partial_z$ near the boundary $\partial (D\times I)$.
\item[(ii)] There is a trajectory of $X$ that enters the plug at
$D\times\{-\varepsilon\}$ and is trapped, i.e.\ it never leaves the plug.
\item[(iii)] The flow of $X$ on $D\times I$ is aperiodic, i.e.\ it does
not have any closed orbits.
\item[(iv)] Any orbit that traverses the plug enters and exits the
plug at a pair of matching points $(x,\pm\varepsilon)$.
\end{itemize}

Inside the plug, there will be an aperiodic invariant set that serves
as a \emph{trap} for at least one orbit; this orbit enters the plug
and becomes asymptotic to the invariant set. By inserting the plug into
a flow box around a point on an isolated periodic orbit
of a given flow such that the trapped trajectory matches
the entrance point of the periodic orbit, one destroys this
periodic orbit without creating any new closed trajectories.

In general, $D$ may be replaced by any $(m-1)$-dimensional
manifold with boundary such that $D\times I$ embeds into $\R^m$,
with $\{p\}\times I$ mapping to a line segment
parallel to the $x_m$-direction for all $p\in D$.
For instance, in the case $m=3$ one can take
$D$ to be any orientable surface with boundary.

There are two features of our plug that make it considerably simpler
than those used by Ginzburg~\cite{ginz95,ginz97} and
Kerman~\cite{kerm02}.
First of all, we use an irrational flow on a torus as
a trap, as in Wilson's original construction, whereas Ginzburg
worked with the horocycle flow on the unit tangent bundle
of a hyperbolic surface, and Kerman with the dense aperiodic subset
inside the  geodesic flow on the unit tangent bundle of a torus.

Secondly, the main task for both Ginzburg and Kerman is to construct a
symplectic embedding of the plug into a flow box of the
original Hamiltonian flow. For Ginzburg's plug, this
requires a subtle application of Gromov's $h$-principle,
see also~\cite[pp.~118--120]{elmi02}. In our
construction, the plug comes from
a deformation of a contact hypersurface in standard symplectic
space, so the symplectic embedding comes for free.

A flow on $D\times I$ satisfying properties (i) to (iii) will be
called a \emph{half-plug}. Our construction in \cite{grz14} yields
a half-plug for Reeb flows. Here we show that by taking
its mirror image under $z\mapsto -z$ and reversing the flow direction,
this is still a half-plug for Hamiltonian flows. When one
half-plug is put on top of the other, the matching condition (iv)
will be satisfied. This is reminiscent of Wilson's construction.
\section{The trap for Reeb orbits}
A \emph{contact form} on a $(2n-1)$-dimensional manifold is
a $1$-form $\alpha$ such that $\alpha\wedge(\rmd\alpha)^{n-1}$
is a volume form. An example is the standard contact
form
\[ \alst=\rmd z+\frac{1}{2}\sum_{j=1}^{n-1}
(x_j\,\rmd y_j-y_j\,\rmd x_j)\]
on $\R^{2n-1}$.
The \emph{Reeb vector field} of  a contact form $\alpha$
is the unique vector field $R$ satisfying
$\rmd\alpha(R,\,.\,)\equiv 0$ and
$\alpha(R)\equiv 1$.
The Reeb vector field of $\alst$ is $\partial_z$.

For $n\geq 3$, in \cite{grz14} we constructed a trap for Reeb orbits as the
Reeb flow of a suitable
contact form $\alst/H$, where $H\co\R^{2n-1}\rightarrow\R^+$
is a smooth function that is identically $1$ outside a compact set.
This function may be chosen $C^0$-close to~$1$. The map
\[ (x_1,y_1,\ldots,x_{n-1},y_{n-1},z)\longmapsto (\lambda x_1,\lambda y_1,
\ldots,\lambda x_{n-1},\lambda y_{n-1},\lambda^2 z),\]
with $\lambda\in\R^+$, pulls back $\alst$ to $\lambda^2\alst$.
This rescaling allows one to choose the support of $H-1$ in
an arbitrarily small neighbourhood of the origin.

We showed that a function
$H$ can be found such that the Reeb flow of $\alst/H$ remains aperiodic,
and some Reeb orbits become asymptotic to
an irrational flow on a Clifford $(n-1)$-torus, and hence trapped.
So this local model is a half-plug,
but the matching condition~(iv) is not satisfied.
Indeed, as there are cases where the Weinstein conjecture for Reeb flows
has been resolved positively, there can be no general
plug construction for Reeb orbits.
\section{Hamiltonian interpretation of the half-plug}
Now let $(W,\omega)$ be a $2n$-dimensional symplectic manifold,
i.e.\ $\omega$ is a closed $2$-form on $W$ such that $\omega^n$
is a volume form. Let $K\co W\rightarrow\R$ be a smooth function.
Then the \emph{Hamiltonian vector field} $X_K$ corresponding to $K$
is defined by
\[ \omega(X_K,\,.\,)=-\rmd K.\]
If $\Sigma=K^{-1}(c)$ is a regular level set of~$K$, then $X_K$
is a tangent vector field along this smooth hypersurface.
By replacing $K$ by $K-c$, we may always assume that $\Sigma$
is the zero level set of~$K$. Given any other function $K'$
with $\Sigma$ as its regular zero level set, we have $K'=fK$ with
$f$ some smooth nowhere zero function, and then $X_{K'}=fX_K$.
Hence, up to reparametrisation, the Hamiltonian flow is determined
by~$\Sigma$ and $\omega$ alone.

In fact, a little more is true.
The $2$-form $\omega$ restricts to a $2$-form $\omega_{\Sigma}:=
\omega|_{T\Sigma}$ of maximal rank $2n-2$, and the flow lines of $X_K$
are the characteristics of $\omega_{\Sigma}$, i.e.\ trajectories tangent
to the kernel of~$\omega_{\Sigma}$. By the symplectic neighbourhood
theorem~\cite[p.~104]{mcsa98}, a neighbourhood of any
oriented hypersurface $\Sigma$ in $(W,\omega)$ is symplectomorphic to 
$(-\sigma,\sigma)\times\Sigma$ with symplectic form
$\omega_{\Sigma}+d(s\beta)$, where $s$ is the coordinate in
$(-\sigma,\sigma)$, and $\beta$ is a $1$-form on $\Sigma$
that does not vanish in the characteristic direction. Then the vector
field $X$ in this characteristic direction (i.e.\ in the
kernel of~$\omega_{\Sigma}$) with $\beta(X)=1$ is the
Hamiltonian vector field corresponding to the function~$s$.
In other words, the Hamiltonian flow is completely determined
(up to reparametrisation) by $(\Sigma,\omega_{\Sigma})$.

The \emph{symplectisation} of a contact manifold $(M,\alpha)$
is the symplectic manifold $\bigl(\R\times M,\rmd(\rme^t\alpha)\bigr)$.
The Hamiltonian vector field corresponding to the function~$\rme^t$
at the level $\rme^0=1$
is then the Reeb vector field $R$ of~$\alpha$. The rescaled contact form
$\rme^f\alpha$ on $M$, where $f$ is some smooth function on~$M$,
can be obtained by pulling back the $1$-form
$\rme^t\alpha$ under the embedding $M\ni x\mapsto(f(x),x)\in\R\times M$.
The Reeb vector field of $\rme^f\alpha$, when interpreted as
a vector field along that graph embedding, is the Hamiltonian vector field
of the function $\rme^{t-f}$ at the level~$1$.

After these preliminaries, we now want to interpret our Reeb trap
in this Hamiltonian setting. Consider $\R^{2n}$ with the standard
symplectic form
\[ \omst=\rmd w\wedge\rmd z+\sum_{j=1}^{n-1}\rmd x_j\wedge\rmd y_j.\]
The vector field
\[ Y= \frac{1}{2}(w\partial_w+z\partial_z)+\frac{1}{2}\sum_{j=1}^{n-1}
(x_j\partial_{x_j}+y_j\partial_{y_j})\]
is a Liouville vector field for $\omst$, i.e.\ $L_Y\omst=\omst$,
and $\R^{2n-1}\equiv\{w=2\}$ is transverse to $Y$
and hence a contact type hypersurface, on
which $i_Y\omst$ restricts to~$\alst$.

The symplectisation $\bigl(\R\times\R^{2n-1},\rmd(\rme^t\alst)\bigr)$
of $(\R^{2n-1},\alst)$ can be identified symplectically with
$(\{w>0\},\omst)$ by identifying $\{0\}\times\R^{2n-1}$
with $\{w=2\}$ and mapping the flow lines of $\partial_t$ to those of~$Y$.
Now, replacing $\alst$ on $\R^{2n-1}$ by $\alst/H$, as in the construction
of our trap, amounts to replacing $\{0\}\times\R^{2n-1}$ by the
graph of $-\log H$ in the symplectisation, and the new Reeb flow corresponds
to the Hamiltonian flow on this new embedding of $\R^{2n-1}$. This
embedding is isotopic (under a compactly supported
isotopy) and $C^0$-close to the original one;
the analogous statement holds for the corresponding embeddings in
$(\R^{2n},\omst)$.

Now suppose we have an isolated periodic orbit $\Gamma$ in the Hamiltonian
flow of $X_K$ on $\Sigma=K^{-1}(c)$. Since $\omega_{\Sigma}$
is invariant under the flow of $X_K$, Darboux's theorem allows us to choose a
flow box $B=D^{2n-2}\times I$ around a point on $\Gamma$ such
that $\omega_{\Sigma}$ is given by
\[ \sum_{j=1}^{n-1}\rmd x_j\wedge\rmd y_j=\rmd\alst\]
on the flow box. The symplectic form in a neighbourhood
$(-\sigma,\sigma)\times B$ is then given by
\[ \omega_{\Sigma}+\rmd(s\alst)=\rmd s\wedge\alst+(1+s)\,\rmd\alst.\]
With the substitution $s=\rme^t-1$, this becomes the symplectic form
on the symplectisation of $(B,\alst)$.

By the construction of the Reeb trap, we can trap the periodic orbit
$\Gamma$ by a deformation of $\{0\}\times B$
inside $(-\sigma,\sigma)\times B$, supported in the interior of~$B$.
\section{A Hamiltonian plug}
In order to build a Hamiltonian plug, we also need to take care
of the matching condition~(iv). To this end, we think of
two boxes $B_{\pm}$ sitting inside~$B$. With $B=D^{2n-2}_{\delta}\times
[-\varepsilon,\varepsilon]$, we take
\[ B_+=D^{2n-2}_{\delta/2}\times[-3\varepsilon/4,-\varepsilon/4]\]
and
\[ B_-=D^{2n-2}_{\delta/2}\times[\varepsilon/4,3\varepsilon/4],\]
say. On $B_+$ we perform the previous construction, so we replace
the linear flow in the $z$-direction by the flow of
the Reeb vector field $R$ of $\alst/H$,
realised as a Hamiltonian flow by a deformation of $\{0\}\times B_+$
inside $(-\sigma,\sigma)\times B_+$, supported in the interior of~$B_+$.

Condition (iv) will be satisfied if the linear flow on $B_-$
is replaced by the flow of $-\Phi^*R$, where $\Phi(\bfx,\bfy,z)=
(\bfx,\bfy,-z)$, i.e.\ the negative
Reeb flow of the contact form $\Phi^*(\alst/H)$.
In other words, we simply reverse the Reeb flow
in our local model, and turn the local model upside down.
Notice that $\Phi^*\rmd\alst=\rmd\alst$, 
so, by the neighbourhood theorem for hypersurfaces,
$\Phi$ extends to a symplectomorphism of a neighbourhood of 
$B_-$ to a neighbourhood of~$B_+$ in $W$. This diffeomorphism,
however, is not straightforward, since $\Phi^*$ does not pull back
$\alst$ to $-\alst$. It is clear, however, that this
extended symplectomorphism must reverse the coorientation of $B$ in~$W$.

The negative Reeb flow, too, can be interpreted as
a Hamiltonian flow; in the symplectisation it would be
the one corresponding to the Hamiltonian function $-\rme^t$ at the
level~$1$. So the desired flow on $B_-$ can likewise be realised
by a deformation of the hypersurface in the symplectic manifold.

Here , briefly, is an alternative look at this mirror construction.
By the symplectic neighbourhood theorem
we may write the symplectic form on a neighbourhood
$(-\sigma,\sigma)\times B_+$ of $B_+$ in $W$ as
$\rmd s\wedge\rmd z+\rmd\alst$; a perhaps smaller
neighbourhood is symplectomorphic to
a neighbourhood of $(\{0\}\times B_+,\alst)$ in its symplectisation
$\bigl(\R\times B_+,\rmd (\rme^t\alst)\bigr)$.
In this neighbourhood we perform the deformation to produce a half-plug.

The symplectic form on a neighbourhood $(-\sigma,\sigma)\times B_-$
of $B_-$ in $W$ can likewise be written as $\rmd s\wedge\rmd z+\rmd\alst$.
This contains a smaller neighbourhood symplectomorphic to
the symplectisation of $(\{0\}\times B_-,\Phi^*\alst)$, where
positive values of $t$ in the symplectisation correspond to negative
values of~$s$. The negative Reeb field of $\Phi^*\alst$ is $\partial_z$,
and after the `mirror' deformation, the negative Reeb flow will be as desired.
\section{Aperiodic Hamiltonian and volume-preserving flows}
One can now, as in \cite{ginz97} and \cite{kerm02},
construct smooth aperiodic Hamiltonian flows on hypersurfaces
of dimension at least five by starting with a Hamiltonian flow
having only isolated periodic orbits. Examples are ellipsoids
$\bigl\{\sum_{j=1}^n a_j|z_j|^2=1\bigr\}$ in $\C^n$ with
$a_1,\ldots, a_n$ positive rationally independent real numbers.
Non-simply connected hypersurfaces with this property
have been constructed by F.~Laudenbach~\cite{laud97}.

The flow of a Hamiltonian vector field $X_K$ on a 
$(2n-1)$-dimensional hypersurface
$\Sigma=K^{-1}(c)$ preserves the volume form
$\beta\wedge\omega_{\Sigma}^{n-1}$, where $\beta$ is
a $1$-form on $\Sigma$ with $\beta(X_K)=1$. More generally,
our plug can be inserted into a flow box of a volume-preserving
flow on any manifold of dimension $2n-1\geq 5$, with the volume form
on the plug defined by the contact form.
In even dimensions $2n\geq 6$, we take the product of
our half-plugs in dimension $2n-1$ with an interval
$[-\varepsilon,\varepsilon]$. On this interval
we consider a smooth function $\psi\co
[-\varepsilon,\varepsilon]\rightarrow [0,1]$, supported in
$[-\varepsilon/2,\varepsilon/2]$ and with $\psi(u)=1$ for $u$ near zero.
Then, on the slice $B_+\times\{u\}$, we take the
Reeb flow of the contact form
\[ \alpha_u:=\alst/\bigl(\psi(u)\cdot H+1-\psi(u)\bigr).\]
This flow preserves the volume form $\Omega$ on
$B_+\times[-\varepsilon,\varepsilon]$ given 
along the slice $B_+\times\{u\}$ by
$\alpha_u\wedge(\rmd\alpha_u)^{n-1}\wedge\rmd u$. On $B_-\times
[-\varepsilon,\varepsilon]$ we mirror this construction.

By \cite[eqn.~(H-iv)]{grz14}, the Reeb vector field $R_u$ of $\alpha_u$
satisfies
\[ \rmd z(R_u)=\psi(u)H-\frac{\psi(u)}{2}\sum_{j=1}^{n-1}
\bigl(x_j H_{x_j}+y_jH_{y_j}\bigr) +1-\psi(u),\]
and $H$ was chosen such that for $\psi(u)=1$ this expression is
zero on the Clifford torus and positive elsewhere. Thus, for $\psi(u)\in [0,1)$,
we have $\rmd z(R_u)>0$. This guarantees that we do not produce
any new periodic orbits inside the half-plug.

Thus, again, starting from any volume-preserving flow in
dimension at least five with isolated periodic orbits, one can produce
a smooth aperiodic one. According to~\cite{ginz99},
Wilson's plug can be chosen divergence-free in
all dimensions $\geq 4$. On $3$-manifolds,
aperiodic volume-preserving flows of class $C^1$
have been constructed by G.~Kuper\-berg~\cite{kupe-g96}.
\begin{ack}
We thank the referee of our earlier paper for valuable
suggestions,
which prompted the present note. We also thank Peter Albers
for useful conversations.
H.~G.\ and K.~Z.\ are partially supported by DFG grants GE 1245/2-1
and ZE 992/1-1, respectively. N.~R.\ is supported by post-doctoral grants
from the Mittag-Leffler Institute and SFB 878 `Groups, Geometry
and Actions' at WWU M\"unster.
\end{ack}

\end{document}